\input amstex

\documentstyle{amsppt}

\magnification = 1200

\baselineskip=24pt  

\define\<{<}
\define\>{>}

\define\lra{\longrightarrow}

\define\A{{\Cal A}}

\define\sign#1{\text{sign}\, #1 \,}
\redefine\bar{\widetilde}
\redefine\qed{ \hfill {$\diamondsuit$}}

\document
\redefine\sect{1}

\heading{ $\Cal S$-structures for $k$-linear categories and
the definition of a modular functor}
\endheading
\footnote""{MR-classification: 57N10, 18D10, 81E05, 16A16.}

\hfil {Ulrike Tillmann} \hfil

\vskip .5in

\subheading{ Introduction}

\vskip .1in
\noindent
{\it Motivation and background.}
Motivated by ideas from string theory and quantum field theory
new invariants of knots and 3-dimensional manifolds have been constructed 
from complex algebraic structures such as Hopf algebras  
[17] [22], monoidal categories with additional structure [24],
and modular functors [14]
[23].  These constructions are  closely related. Here we take
a  unifying categorical approach based on  a natural 2-dimensional
generalization of a topological
field theory  in the sense of Atiyah [1], and
show that the 
axioms defining these complex  algebraic structures are a consequence
of the  underlying geometry of surfaces.

\vskip .1in
Recall the following folklore from 1+1 dimensional
topological field theory.
Let $\Cal S$ denote the  oriented
1+1 dimensional cobordism category with monoidal product defined by disjoint 
union. A  topological field theory is
just a monoidal functor from $\Cal S$ into the category of 
finite dimensional vector spaces over the ground field $k$ with tensor
product. In particular, it  
takes the union of two 1-manifolds to the tensor product of  the vector
spaces associated to each  of them. On objects, 
such a functor is then completely
determined by the vector space associated to the unit circle.
It is well known that the $\Cal S$-structure on
this vector space implies that it  is a finite dimensional algebra
with a non-degenerate inner product $<\, , \, >$ such that 
$<ab,c> = <a,bc>$. In other words, it is a Frobenius algebra.
Conversely, any finite dimensional
Frobenius algebra determines a topological field theory.
One gets a topological invariant of 2-manifolds by considering
any closed surface as a morphism from the empty manifold to itself which
induces a linear map from $k$ to itself and hence determines an element in $k$.

\vskip .1in
In order to encode 3-dimensional information, one has to enrich 
$\Cal S$ and make it into a 2-category by adding the group of  orientation 
preserving diffeomorphisms 
of the surfaces as 2-morphsims.  One can now obtain 3-dimensional information
by studying monoidal functors into a higher dimensional analogue of the 
category of finite dimensional vector spaces. This will be a modification
of  the category 
$k$-CAT of $k$-linear categories
in which  each morphism space is 
a vector space over $k$. The category $k$-CAT  is naturally a 
2-category 
where the 1-morphisms are 
linear functors and the 2-morphisms are natural transformations.  
The usual tensor product of vector spaces induces a monoidal structure  on this
2-category.
As before,  on objects
such a monoidal functor will be completely determined  by the category
$\Cal A$ 
associated to the unit circle.
The motivating example is the category of positive energy representations
of the free loop group of a compact Lie group for a given level [18] [19].

\vskip .2in 
\noindent
{\it Results.}
As for TFT's, the action of $\Cal S$ gives $\Cal A$ a rich structure.
We show here that  $\Cal A$ is a semi-simple, Artinian category. Indeed, it is 
equivalent to the category of finitely generated  modules of a semi-simple,
Artinian
algebra. This is a direct consequence of the fact that the cylinder
defines a non-degenerate \lq inner product' on $\Cal A$.
Furthermore, $\Cal A$   
is   balanced,  that  is braided monoidal with 
a compatible twisting. The monoidal structure comes from the pair-of-pants
surface. For the braiding and the twist one considers the appropriate
diffeomorphisms of the pair-of-pants surface and the cylinder.
By analysing the adjoints with respect to the inner product, one
can show that $\Cal A $ has an involution which is compatible with the
twisting.
Similar results hold for any commutative ring $k$.
In analogy to  the  TFT's described above, 
we call such a category Frobenius.
For $k=\Bbb C$, we then analyse the associated Verlinde Algebra and
show that $\Cal A$ gives rise to modular functors, which
according to Walker [23] and Kontsevich [14]
are in one-to-one
correspondence to 2+1 dimensional topological field theories and hence
define 3-manifold invariants.
By extension it seems natural to call 
a monoidal functor from $\Cal S$ to 
$k$-CAT a modular functor.

\vskip .2in
\noindent
{\it Final remarks.}
The class of categories  studied in connection with 3-manifold
invariants has been the class of tortile categories, that is balanced
categories with dual objects (see [20]), a notion that is stronger than that
of Frobenius.
Yetter   [24] constructs 3-manifold invariants from any semi-simple
tortile category thus generalizing the invariants constructed in [22] 
from quantum groups.
The involution in $\Cal A $, however, does  not in general define dual objects,
and therefore $\Cal A$ cannot be tortile.
Our results  suggest that one should study
the notion of a Frobenius category instead.  
Particularly interesting questions are   whether there is
a coherence theorem for Frobenius categories and what extra conditions
a semi-simple, Artinian Frobenius category must satisfy in order to give
rise  to a modular functor. 

\vskip .1in
In this paper we study monoidal functors from $\Cal S$ to $k$-CAT which leads
to a notion of Frobenius categories. 
In principle, one could replace $k$-CAT by any monoidal
category.  
For example, in connection with string theory one
considers the category of cochain complexes  which leads to 
Batalin-Vilkovisky algebras [6].  
On the other hand, monoidal functors to  the stable homotopy category 
give rise to a notion of Floer homotopy type [4].
A most general and rigorous categorical
approach is taken  by Carmody in his recent thesis [3], where one can also find
a proof of the above mentioned TFT-folk theorem. 

\vskip .2in
\noindent
{\it Contents.}
\roster
\item"" 1. Definition of a Modular Functor
\item"" 2. Duality and semi-simplicity of $\Cal A$
\item"" 3. Fusion, involution, and adjoints
\item"{}" 4. The Verlinde algebra 
\item"{}" 5. Orientation and  modular functors
\item"{}" Appendix: Linear Categories
\item"{}" References
\endroster

\vskip .2in
\demo{Acknowlegdements}
This paper started out as a joint project with Graeme Segal.
I would like to thank him for his generosity with his ideas and time.
I am  also grateful  to  Martin
Hyland for many fruitful discussions and  tutorials
on category theory.
Finally, I would like to thank Jonathan
Morris for producing photo-ready drawings
for all figures.
\enddemo

\vskip .4in
\redefine\sect{1}
\redefine\A{ k \text{-CAT}}
\define\wA{\widehat \A}

\subheading {\sect . Definition of a  Modular Functor}

\vskip .1in
We first define the cobordism 2-category $\Cal S$
and then define a modular functor  $\Cal A $
as a monoidal functor from  $\Cal S$
to the  2-category $\wA$
of $k$-linear, amenable categories, $k$-linear functors and
natural transformations.
As the functor $\Cal A$ is monoidal, it is determined on objects by its
image of the unit circle $ S ^1$. Thus a modular functor determines
a linear  category  with an action of the 1+1 dimensional 
cobordism category $\Cal S$.  

\vskip .2in
\subheading { 1.1. The 1+1 dimensional cobordism 2-category}
Heuristically, the 1+1 dimensional cobordism category  is the category
of all closed, oriented, compact  1-manifolds and  oriented, compact
2-manifolds and their diffeomorphisms;
disjoint union of manifolds induces a symmetric monoidal structure.
Unfortunately, so far there is no satisfactory categorical concept
to capture all its features.
\footnote{
An attempt to describe the most general cobordism category may be found in
[3].}
In order to stay \lq honest', we will ignore diffeomorphisms of 1-manifolds
(but compare Section 5) and only take into account the connected components
of the diffeomorphisms of 2-manifolds. More precisely, we will work with
the model $[\Cal S _\Gamma]$ constructed in [21] which is a monoidal 
strict  2-category 
in the sense of [7] and has an obvious symmetry.
\footnote{
$[\Cal S _\Gamma ]$ has a braiding in the sense of [13] and is symmetric in the
sense of [3]. However, as long as there are no coherence theorems for 
symmetric or braided monoidal 2-categories there is no need to be
precise about this here. 
On the other hand, $[\Cal S _\Gamma]$ is so close to being
symmetric strict monoidal (the subcategory $[\Cal S _\Gamma ^b]$ without closed
surfaces is!) that any good notion of a symmetric monoidal
2-category will include $[\Cal S _\Gamma ]$. The mere existence of this model
also suggests that  the abstract 1+1 dimensional cobordism
category is equivalent to $[\Cal S _\Gamma ]$ in an appropiate sense so that
there is not much loss of generality when considering this model.}
Thus define $\Cal S := [\Cal S _\Gamma ]$ and recall its definition.

\vskip .1in
The objects of $\Cal S$ are in one-to-one correspondence with the natural 
numbers, where $n\in \Bbb N \cup \{ 0 \}$ represents 
$n$ circles of radius $1/4$ and
centered at $(1,0), \dots , (n,0)$ in the $xy$-plane in $\Bbb R^3$, or the
empty set.  

A 1-morphism from $m$ to $n$ is a compact,  smooth surface $\Sigma
\subset \Bbb R^2 \times [0,t]$ with $m$ and $n$ boundary circles in the planes
$\Bbb R ^2 \times \{ 0\}$ and $\Bbb R^2 \times \{ t\}$ respectively.
Composition  $\Sigma _2 \circ \Sigma _1$
is defined by placing $\Sigma _2$ on top
of $\Sigma _1$ (by using rigid shifts along the $z$-axis). All
our surfaces are assumed to have collars so that composition  again yields
a smooth surface. 
Straight cylinders 
represent identities after certain identifications have been
made
(see below).  Note that any 1-morphism has a canonical orientation so that
its normal vector is outward pointing.

2-morphisms between $\Sigma _1 $ and $\Sigma _2$ are the connected components
$\Gamma (\Sigma _1 , \Sigma _2) := \pi _0 \text {Diff }^+ ( \Sigma _1,
\Sigma _2; \partial)$ of orientation preserving diffeomorphisms that fix
a collar  (a neighbourhood of the  
boundary of $\Sigma _1 $ and $\Sigma _2$ suitably identified via a rigid
shift along the $z$-axis) pointwise.
The set of 2-morphisms have two compositions: an internal composition
which is induced by composition of diffeomorphisms, and an external
composition which is induced by the composition of 1-morphisms. As the 
diffeomorphisms fix the boundary pointwise, this latter composition is also
well-defined.

A monoidal structure  $\otimes $ is induced by disjoint union.
On objects this is $n \otimes m := n+m$. In order to define $\Sigma _1 \otimes
\Sigma _2 $, first the height of $\Sigma _2$ must be adjusted to the height
of $\Sigma _1$ by a chosen
reparametrization of its $z$-coordinates. Then $\Sigma _2$
needs to be  shifted  to the right
along the $x$-axis by possibly different amounts for each
$z$-coordinate.

It is clear that $\otimes $ is neither associative nor functorial.
To solve this problem we need to identify all surfaces which can be mapped 
to each other by finite repetitions of the following procedure:
Cut  a surface $\Sigma  $ by a finite number of horizontal planes into connected
components $\Sigma _ i$; Each $\Sigma _i$
may be deformed independently by  a reparametrization of the
$z$-coordinate or by  a smooth \lq shift'
along the $x$-axis by different amounts for each $z$-coordinate but
fixing the boundary of $\Sigma _i$.

There are two important features of this identification. First of all,
note that  different components between two horizontal
planes can be moved through each other. Thus twisted cylinders define
not just a braiding but also
a symmetry for the monoidal structure. Secondly, any two surfaces so
identified are not just diffeomorphic but the procedure above defines a 
unique diffeomorphism up to isotopy.

\vskip .2in
\proclaim{Theorem \sect .1 [21]}
$\Cal S$ is a semistrict monoidal strict 2-category with a symmetric braiding.
\endproclaim

\vskip .1in
\demo{Ad figura} 
Occasionally we will use diagrams to indicate natural diffeomorphisms
between two surfaces $\Sigma _1 $ and $\Sigma _2$ in $\Bbb R ^3$. 
When these are genus
0 surfaces, as they  will be most of the time, 
any two isotopies (relative to $\Bbb R ^3$)
differ at most by an isotopy of the target surface (relative to itself),
and hence define a unique  2-morphism.
However, when the genus is not zero this is no longer true!
We trust that also in these cases the picture will unambiguously
suggest  the  
definition of a 2-morphism.
\enddemo

\vskip .2in
\subheading {1.2. Definition of a modular functor}
Given a commutative ring $k$,
let $\wA $ be the  strict 2-category of amenable $k$-linear categories.
Its objects are additive and idempotent complete
$k$-linear categories, 
its 1-morphisms are $k$-linear functors, and its 2-morphisms are 
natural transformations.
When $k$ is a field,  $\wA$ is also called a 2-vector space.
The usual tensor product over $k$ induces a symmetric monoidal structure
on $\wA$.
\footnote{
Similar remarks as in the previous footnote also apply to $\wA$.}
The definition of amenable and the construction of the tensor product
as well as further details may be found in the Appendix.

\vskip .1in
Basically, a modular functor is a monoidal functor $\Cal A :\Cal S \to \wA$.
Because of the complicated structure of the underlying categories, there are
many possible choices depending on how much structure is to be preserved 
on the nose, up to isomorphisms, up to natural transformations, $\dots$.
We will assume here that $\Cal A $ is a homomorphism of 2-categories 
which preserves the monoidal structure up to pseudo-natural equivalences
(see [12]).
\footnote{
This is a special case of a trihomomorphism in the sense of [7]. We will not 
need the full strength of this definition.}
Thus, to each object
$n\in \Cal S$, it assigns an amenable $k$-category $\Cal A (n)$, to a cobordisms
$\Sigma \in \Cal S (n,m)$ a functor
\footnote { We will see later that there is a natural equivalence of functor
categories
$$
\Cal A (n) ^{\Cal A(m)} \simeq (\hat k ^{\Cal A (n) ^{op}} ) ^{\Cal A (m)} 
\simeq \hat k ^{ \Cal A (m) \otimes \Cal A (n)^ {op} }
$$
Hence $F_\Sigma \in \Cal A (n) ^{\Cal A  (m)}$  naturally gives rise to an
$\Cal A (n) \otimes \Cal  A (m) ^{op}$-module.  This is the point of view taken
in [2].}
$$
F_\Sigma := \Cal A (\Sigma ) : \Cal A (m) \lra \Cal A (n)
$$
and to each diffeomorphism $\phi : \Sigma _1 \to \Sigma _2 $ a natural
transformation $\Cal A (\phi)$.
Furthermore, the monoidal structure is preserved. In 
particular  there
are  isomorphisms of categories
$\Cal A (n+m) \simeq \Cal A (n) \hat \otimes \Cal A (m)$ which are well-behaved
with respect to 1- and 2-morphisms.
Finally, we will also assume that under these isomorphisms the image 
$\Cal A (T)$
of a pair of twisted cylinders  corresponds
to the  twist functor $\tau: \Cal A (n) \hat \otimes \Cal A (m) \simeq
\Cal A (m)
\hat \otimes \Cal A (m)$.

\vskip .1in
In what follows we will suppress all structure maps and
natural isomorphisms. Also without  loss of generality, we may assume that
the monoidal unit in $\Cal S$ is mapped to the monoidal unit in $\wA$.
Hence,  we simply write $\Cal A (n) \simeq \Cal A (1) ^{\hat \otimes  n}$ and
$\Cal A (\emptyset) = \hat k$. For a cobordism from
$0$ to $n$ there is a well-defined object
$$
E_\Sigma := F _\Sigma (k) \quad \in \Cal A (n)
$$
where $k$ is the canonical object in $\hat k$.

\vskip .1in
Let $I \in \Cal S (1, 1)$ be the standard torus of length 
one, and $F_I:\Cal A(1) \to \Cal A(1)$ 
be the corresponding
functor. Clearly, $F_I \circ F_I \simeq F_I$ as $ I \circ I = I $.

\vskip .2in
\proclaim {\sect .2 Lemma} 
Up to  natural isomorphism, 
$F_I$ is a projector on $\Cal A$.
\hfill { \qed}
\endproclaim
In particular $F_I$ is  naturally equivalent to the identity 
on the subcategory generated
by the elements $E_\Sigma$ associated to surfaces $\Sigma$ as
$$
F_I (E_\Sigma) = F_I \circ F_\Sigma (k) \simeq F_\Sigma (k).
$$
As we are mainly interested in analysing 
the structure of $\Cal A (1)$ induced by the
$\Cal S$-action coming from surfaces, we will restrict our attention
to the image of $\Cal A (1)$ under the projection  $F_I$. 
This leads us to the following definition.

\vskip .2in
\proclaim{\sect .3 Definition}
A modular $k$-functor is a monoidal functor  $\Cal A :\Cal S \to 
\wA $ in the above sense
such that $F_I : \Cal A (1) \to \Cal A (1) $ is 
naturally equivalent to
the identity functor. We say, the $k$-linear category
$\Cal A  := \Cal A (1)$ has an $\Cal S$-structure.
\endproclaim

\vskip .4in
\redefine\sect{2}

\subheading { \sect . Duality and semi-simplicity of $\Cal A$}

\vskip .1in
We  prove that when $k$ is a field the category $\Cal A= \Cal A (1)$
is semi-simple for any modular functor.
More precisely, we  show that $\Cal A$ is equivalent to the category of
finitely generated projective modules over a separable $k$-algebra  which is 
semi-simple when 
$k$ is a field. 
The geometry of a modular functor is
here used to establish an analogue of
a non-degenerate inner product on $\Cal A$,
and hence self-duality.  Semi-simplicity or
separability of the algebra in question now follow from results in algebra.

\vskip .2in
\subheading{\sect .1 Duality of $\Cal A $} 
Let $C \in \Cal  S (0,2)$ denote a fixed, bent  cylinder 
with two ingoing boundary components.
We will denote the induced functor $F_C : \Cal A \hat \otimes \Cal A
\to \hat k$ also by $< \, , \, >$, thinking of it as an inner product on $
\Cal A$.
The cylinder $\bar C$, $C$ rotated by $180^o$,
defines a map  in the opposite direction.  
Let $E_{\bar C} \in \Cal A \hat
\otimes  \Cal A$ 
be the image of the canonical element $k$ in the category $\hat k$ 
and  
write it as
$$
E_{\bar C} =F_{\bar C} (k)  = (\Sigma _{i=1} ^n ( P_i \otimes Q_i ); \pi)
$$
where $\pi$ is a projector and $Q_i$ and $P_i$ are objects in $\Cal A$.
In this section, a special role is played by  the element
$X:= \Sigma _{i=1}^n Q_i \in \Cal A $ and
its algebra of morphims $A := \Cal A (X,X)$.

\vskip .1in
Let  $\Cal A ^* = (\hat k) ^{\Cal A} $ be the category of $k$-functors
from $\Cal A $ to  the category $\hat k$  of
finitely generated projective $k$-modules.
Define two covariant  functors:
$$
\aligned
\Cal I: \Cal A \to \Cal A ^*  \, &\text{ via } \,  X \mapsto F_C(X,\, ) 
		= 	<X, \, > \\
\Cal E : \Cal A ^* \to \Cal A \, &\text { via } \,  
		F \mapsto (F \hat \otimes id _{\Cal A})  (E_C)
		= ( \Sigma _{i=1} ^n  (F (P_i) \otimes Q _i);  
			(F\otimes id_{\Cal 
					A} ) (\pi)).
\endaligned
$$

\vskip .2in
\proclaim{Theorem \sect .1}
The functors $\Cal I$ and $\Cal E$ are inverse equivalences of categories.
In particular, $\Cal A \simeq \Cal A ^*$.
\endproclaim

\demo {Proof}
This  follows immediately from the assumption that the functor associated to
the cylinder $F _I : \Cal A \to \Cal A $ is an equivalence of $\Cal A$.
More precisely, there is the following string of natural isomorphisms.
(See Figure \sect .1)
$$
\aligned
\Cal E \circ \Cal I (\,) &= \Cal E ( F_C (\, ,\_  ))	\\
&=  (F_C  \otimes id_{\Cal A}) ( \, , E_{\bar C})	\\
&= (F_C  \otimes id _\Cal A) (id _\Cal A\otimes F_{\bar C}) ( \, , k)	\\
&\simeq (F_C \otimes F_I )(F_I \otimes F_{\bar C}) ( \, ) \\
&\simeq F _I ( \, )	\\
&\simeq id_\Cal A (\, )
\endaligned
$$
We have only used the functoriality of the modular functor and the natural 
equi\-valence $\Cal A \hat \otimes \hat k \simeq \Cal A$.
Similarly, we have the following isomorphisms of functors.
$$
\aligned
\Cal I \circ \Cal E ( \, ) &\simeq \Cal I ( ( \,  \hat \otimes id_{\Cal A}) 
			(F_{\bar C} k) ) \\
&= (\, \otimes F_C ) (F _{\bar C} \otimes  id_\Cal A) \\
&= (\, ) (id_\Cal A  \otimes F_C) (F _{\bar C} \otimes id _\Cal A) \\
&\simeq (\, ) F_I = F _I ^* (\, ) \\
&\simeq id _{\Cal A ^*} (\, )
\endaligned
$$
as all functors are linear and $F^*_I \simeq id_{\Cal A ^*}$ by adjointness
(see A.2).
\qed
\enddemo 

\vskip 2in
\hfil{\it Figure \sect .1: $\Cal E \circ \Cal I \simeq F_I \simeq id_{\Cal A}
$.}
\hfil

\vskip .2in
\subheading{\sect .2 Semi-simplicity when $k$ is a field}
Our next goal is to prove  Theorem \sect .5. As preparation we derive the 
following  three Lemmata from Theorem \sect .1.
Recall that $X = \Sigma Q_i$ and $A = \Cal A (X, X )$.

\vskip .1in
\proclaim {Lemma \sect .2 } 
$\Cal A$ is equivalent to the category of finitely
generated projective $A$-modules.
\endproclaim

\demo{Proof}
The category of finitely generated projective
$A$-modules is the
additive and idempotent completion of $\Cal A (X,X)$ considered as
a category with one object $X$, and is therefore naturally a
subcategory of $\Cal A$ as $\Cal  A $ is amenable. 
On the other hand,  for every objet $Y \in \Cal A $, 
$$
Y\simeq \Cal E \circ \Cal I (Y) = ( \Sigma_i <Y,P_i> \otimes Q_i; <id_Y, \pi>)
$$
where $<id_Y, \pi>$ is shorthand for $(<\, , \, >\hat \otimes id_{\Cal A} )
( id _Y \otimes \pi)$. Hence,  interpreting $<Y, P_i> $ as the multiplicity 
of $Q_i$,
$Y$ is isomorphic to a submodule of direct sums of $X$.
Note that
$<\, , \, >$ takes values in $\hat k$, the category of finitely generated
projective $k$-modules. 
Thus,  $\Cal E \circ \Cal I ( Y)$ is a direct summand of finitely
many copies of $X= \Sigma _{i=1} ^n Q_i$. In other words, it
lies in the subcategory
of finitely generated projective $A$-modules.
By the above theorem, $\Cal E \circ \Cal I$ is an equivalence, and therefore
$\Cal A$ is equivalent to the category of finitely generated projective
$A$-modules.
\qed
\enddemo

\vskip .2in
\proclaim {Lemma \sect .3 }
For all objects $Y$ and $Z$
in $\Cal A$, the morphism set  $\Cal A (Y, Z)$ 
is a finitely generated $k$-module.
In particular, $A$ is finitely generated over $k$.
\endproclaim

\demo{Proof}
Under the equivalence $\Cal E \circ \Cal I$, 
$\Cal A (Y, Z)$ is isomorphic to
$$
\{ g \in \Sigma _{i,j} \hom (<Y, P_i>, <Z,P_j>) \otimes \Cal A (Q_i, Q_j)\, |
	\, g<id_Y, \pi> = <id_Z, \pi>g\}
$$
The isomorphism is explicitly given by
$$
f\lra \tilde f <id_Y, \pi> = <id_Z,\pi> \tilde f 
$$
where $\tilde f = \Sigma _i <f, id_{P_i}> \otimes id_{Q_i}$. As each of the
$<Y, P_i>$ and $<Z, P_i>$ are finitely generated projective  $k$-modules,
the set of $\tilde f$'s also form  a finitely generated $k$-module.
Hence, the isomorphic image of $\Cal A (Y,Z)$ is a finitely generated
$k$-module.
\qed
\enddemo

\vskip .2in
\proclaim {Lemma \sect .4 }
If $k$ is a field,
$\Cal A ^* $ is equivalent to  the category of finitely generated
$A$-modules.
\endproclaim

\demo{Proof}
Note that by Lemma \sect .2,
$\Cal A $ is equivalent to the additive 
and idempotent completion of $A$, considered as a category with one element
$X$ and
morphism space $A$.  But every functor $A \to \hat k $ has a unique extension
to $\hat A \simeq \Cal A$, and therefore $\Cal A ^* \simeq \hat k ^A$.

Clearly, a functor $F: A \to \hat k$ is 
determined by the image of $X$ which has to be an $A$-module that is
finitely generated and projective as $k$-module, id est finite dimensional.
Vice-versa, any such module determines a functor $F$.
Similarly,  natural transformations between functors are in one-to-one
correspondence with
$A$-maps of the corresponding  $A$-modules.
But now, since
$A$ itself is finite dimensional by Lemma \sect .3,  
every finitely generated 
$A$-module is finite dimensional over $k$.
\qed
\enddemo

\vskip .2in
\proclaim{Theorem \sect .5 } For  $k$  a field, $\Cal A$ is a semi-simple,
Artinian 
category. More precisely, $\Cal A$
is equivalent  to the category of finitely generated 
modules over the semi-simple, Artinian  algebra $A = \Cal A ( 
X, X )$.
\endproclaim

\demo{Proof} We need to show that $A$ is semi-simple.
Theorem \sect .1 gives us the equivalence $\Cal A \simeq \Cal A ^*$.  By the
above lemmata, $\Cal A$ is equivalent to
the category of projective $A$-modules and
$\Cal A^*$ is equivalent to  the category of finitely generated 
$A$-modules.
Hence, 
every finitely generated $A$-module is projective. 
This can be  true if and only if $A$ is semi-simple.
\qed
\enddemo

\vskip .2in
\subheading {\sect .3 Separability when $k$ is any commutative ring}
In the remainder of this section we will generalize  Theorem 2.5 to
arbitrary commutative rings $k$.  In this case the algebra $A$ turns out to be 
separable. Note that Lemmata 2.2 and 2.3 are valid for any $k$. 

Recall, an algebra $R$ is called separable if it is projective
as $R\otimes R^{op}$-module. This is equivalent to the existence of 
an idempotent $e = \Sigma _i \,  e_i \otimes f_i$ in $R\otimes R^{op}$
such that the multiplication map sends $e$ to $1\in R$ and 
$er=re$ for all $r\in R$. $e$ is called the  separability idempotent. 
We refer the reader to [5] for an extensive study of separable algebras.
From there we quote the following result [5, p.72].

\proclaim {Theorem }
An algebra $R$ over a commutative ring $k$ is separable if and only if
$R_m := R/{mR}$ is a separable $ k_m:=
k/m$-module for every maximal ideal
$m$ of $k$.
\endproclaim

Mitchell  [16] studies a 
natural generalization of the notion of separability to
$k$-categories. It suffices to record here that this notion of
separability 
is invariant under Morita equivalence.

\vskip .2in
\proclaim{ Theorem \sect .6 }
For $k$ a  commutative ring, $\Cal A$ is a separable category. More precisely,
$\Cal A$ is equivalent to the category of finitely 
generated modules of the separable algebra $A = \Cal A (X,X)$ which are
projective as $k$-modules. 
\endproclaim

\demo{Proof} 
Let $\Cal A = \Cal A (1)$  be determined by some $k$-modular functor
$\Cal A: \Cal S \to \wA $.
A  homomorphism $\phi : k \to k'$ of commutative rings
induces a homomorphism of 2-categories
$\Phi = \,\_ \otimes _k k' : \wA \to \widehat{ k' \text{-CAT}}$ which preserves
the tensor product.
Then clearly,
$\Phi \circ \Cal A : \Cal S \to \widehat {k'\text{-CAT}} $ is again a modular 
functor.
Hence, for any maximal ideal $m$ of $k$ and natural projection $k\to k_m$
the ring $A = \Cal A (X,X) $ is taken to $A_m = A \otimes _k k_m
= \Phi \circ \Cal A ( X,X) $. By Theorem \sect .5, $A_m$ is semi-simple
and finite dimensional 
over $k_m$ for all maximal ideals $m$ of $k$. By a similar argument also
$A \otimes_k k'_m $ is semi-simple for all field extensions $k'_m$ of
$k_m$. Hence,  $A_m$ is separable [16, p.49], and 
by the above quoted  theorem,  $A$ is a separable $k$-module.

Now, by Lemma \sect .2, $ \Cal A$ is equivalent to the category of finitely
generated projective $A$-modules. But for a separable algebra $R$, every
module that is projective as $k$-module is also projective as $R$-module
(see [5; p.48]).  Hence the statement of the theorem.  Finally, 
$\Cal A$ is the idempotent completion of $A$ considered as a category,
and thus, by Morita equivalence,  is also separable in the sense of [16].
\qed
\enddemo

\vskip .2in
We end this section by listing some well-known examples 
of separable algebras:

1. Let $R$ be the group algebra $k[G]$ of a finite group $G$ such that
the order of the group $|G|$ is invertible in $k$.  The element $e= \frac 1 
{|G|} \Sigma _{g\in G} \,  (g \otimes g^{-1})$ is a separability idempotent.

2. The algebra of $n\times n$-matrices over $k$.  The element $e = \Sigma _{i,j}
\, e_{ij} \otimes e_{ji}$ is a separability idempotent where $e_{ij}$ is the
$(i,j)$-th  elementary matrix.

3. Similarly,  $R= k ^n$ with componentwise addition and multiplication is
separable. A separability idempotent in this case is given by
$e= \Sigma _i \,  e_i \otimes e_i$ where $e_i $ is the $i$-th generator. 
The algebroids that arise in [2] are all Morita equivalent to 
such a direct product algebra over the ground ring where $n$ is related 
to the number of Jones-Wenzel idempotents of the theory. As separability is 
Morita invariant, this shows that the algebroids themselves are separable.

\vskip .4in
\define\fu{\star}
\define\1{\bold 1}
\redefine\sect{3}

\subheading{\sect . Fusion, involution, and adjoints }

\vskip .1in
$\Cal A$ is shown to be braided monoidal with a compatible twist. 
The duality of the previous section is used to construct an involution
that is compatible with the inner product. Furthermore, relating the adjoint
of a functor to its conjugate under the involution, we are able to show
that the twist is self-dual. Hence, $\Cal A$ satisfies all properties of a 
tortile category except that the image of an object under involution does
not necessarily define a right dual.

\vskip .2in
\subheading { \sect .1 Fusion and balancing in $\Cal A$}
\footnote{
Freed [9] also proves that the category in his somewhat different set-up
is braided monoidal and balanced. The proof here is the same.}
Fix a pair of pants $P$ with two incoming and one outgoing boundary 
components. Also, fix a disk $D$ with one outgoing boundary
circle and define
$$
\aligned
\fu &:= F_P: \Cal A \hat \otimes \Cal A 
\to \Cal A 	\\
\1 &:= E_D \in  \Cal A.
\endaligned
$$ 
A Dehn twist $\phi$ of the cylinder $I$ induces an isomorphism 
of $F_I$, and hence, of  the 
identity functor of $\Cal A$:
$$
\theta :=
\Cal A ( \phi): id_{\Cal A} \to id_{\Cal A}.
$$ 

\vskip .2in
\proclaim{ Proposition \sect .1}
$\Cal A$ is a balanced  category with braided product $\fu$,   unit
$\1$, and twist $\theta$.
\footnote{
These notions have been studied in [11].}
\endproclaim

\demo{Proof}
The functors $\fu (\fu   \otimes id_\Cal A)$ and $\fu (id_\Cal A \otimes \fu)$
are canonically isomorphic to $ F_P (F_P \otimes F_I) $ and $F_P
(F_I \otimes F_P)$. The underlying surfaces are diffeomorphic 
via an isotopy $\phi$. Define the associativity natural isomorphisms
to be $\alpha := \Cal A (\phi)$.
By definition, the two maps of surfaces induced by the pentagon
diagram differ by at most an isotopy relative to the boundary of the
final pairs of pants with four legs. This is the identity
as an element of the
mapping class group and so is its image under $\Cal A$.
Similarly,
we can define a left (right) identity natural isomorphism
$\lambda $ ($\rho$). 
For  $\fu (\1 \otimes id_\Cal A)$ is canonically 
isomorphic to $ F_P (F_D \otimes
F_I)$ which in turn is naturally isomorphic to $F_I$ via an isotopy 
of the underlying surfaces, and hence to  $id _\Cal A$.
Hence, $\fu$ defines a monoidal structure on $\Cal A$.

For the braiding we need a natural isomorphism $\sigma : \fu \to \fu \circ
\tau$ where $\tau: \Cal A \otimes \Cal A \to \Cal A \otimes \Cal A$
interchanges the components. By assumption on the
modular functor, $\tau$ is naturally isomorphic to
$F_T$ where $T$ is some fixed twist of cylinders. Hence we need
a natural isomorphism between $F_P $ and $F_P \circ F_T$ (Figure \sect .1).
Define $\sigma := \Cal A ( \phi )$ where $\phi $ is an isotopy that induces the
obvious half twist necessary.  To check that $\sigma$ is indeed a braiding
we need to check the two hexagon diagrams  for $\sigma$
and $\sigma ^{-1}$.
This corresponds to checking an equation in the mapping class group 
$\Gamma (\Sigma _{0, 3+1})$. 
For $\sigma $ this is done in Figure \sect .2 where we have suppressed the
associativity maps and represented a twist as a permutation of labels.
Compatibility  of $\sigma$ with $\lambda$ and $\rho$, that is 
$\rho \simeq \lambda \sigma$ and $\lambda \simeq
\rho \sigma$,  can be verified in a 
similar way.

To check that $\theta$ is compatible with the braiding   amounts to checking
an equation in the mapping class group $\Gamma (\Sigma _{0, 2 +1})$.
This is done in Figure \sect .3. 
\qed
\enddemo

\newpage
.
\vskip 1.7in
\hfil{\it Figure \sect .1: $\sigma : \fu \overset \sim \to \lra \fu \circ 
\tau $.}\hfil

\vskip 2.5in
\hfil{\it Figure \sect .2: Hexagon relation for $\sigma$.}\hfil

\vskip 2.5in
\hfil{\it Figure \sect .3: $\theta _{A\fu B}= \sigma _{B,A} (\theta _B \fu 
\theta _A) \sigma _{A,B}$ and $ \theta _1 = id _1$.}\hfil

\define\tra{\text {tr}}

\subheading{\sect .2 Involution and compatibility of the inner product}
For a fixed $Y\in \Cal A $ consider the covariant functor $Z\to \Cal A (Y,Z)$.
By the results of Section 2, $\Cal A (Y,Z)$ is a finitely generated
projective $k$-module. Hence, we can define a functor
$$
\hom : \Cal A^{op} \lra  \Cal A ^* \quad \text{via} \quad    Y \mapsto 
\Cal A (Y, \_).
$$
Composing this with the equivalence $\Cal E :\Cal A ^*\to \Cal A$, we get a
contravariant functor of $\Cal A $ to itself:
$$
^* : \Cal A ^{op} \lra \Cal A \quad \text {via} \quad Y\mapsto Y^* := \Cal E
						(\Cal A (Y, \_ )) .
$$
Note that $\hom : \hat k ^{op} \to \hat k ^*$ sends a $k$-module $V$ to 
$\hom _k (V, \_)$ which in turn is mapped by $\Cal E $ to the dual module
$V^* =\hom _k (V, k)$.
\footnote
{ I expect  that in general this involution is equivalent to
the functor that takes an $A ^{op}$-module $M$ to its dual $ \hom _k(M,k)$
which is now an $A$-module.}

\vskip .2in
\proclaim{ Proposition \sect .2} 
$^* : \Cal A ^{op} \to \Cal A $ is an involution on $\Cal A$ in the sense that
its square $^* \circ ^ *$ is equivalent to the identity functor.
\endproclaim

\demo{Proof} 
$\Cal A (Y, \_) : \Cal A ^{op} \to Mod (k)$ is faithful and full.
Thus,  $ A ^{op} = \Cal A ^{op} (X, X)$ is mapped isomorphically to $\Cal A (
X^*,  X ^*)$ where $X = \Sigma _i Q_i$. 
Hence, as $\Cal A $ is equivalent to  the category of
finitely generated $A $-modules, 
it is equivalent  to the category of finitely generated  
$\Cal A (  (X ^*)^*,  (X ^* )^*)$-modules.
More explicitly, the sequence of natural equivalences
$$
\Cal A ^{op} \overset \hom \to \lra \Cal A ^* \overset \Cal E \to \lra
	\Cal A \overset \Cal I \to \lra \Cal A ^* \overset \hom \to 
	\lra (\Cal A ) ^{*^*}
$$
is equivalent to the canonical functor $\Cal A^{op}\to ((\Cal A^{op}) ^*)^*$
via $\Cal I \circ \Cal E \simeq id_{\Cal A^* }$.
\qed
\enddemo

\vskip .2in
\proclaim{Proposition \sect .3} 
There is a natural equivalence 
$ 
\Cal A (\_ ,\_) \simeq < \_ ^*, \_ > .
$
\endproclaim

\demo {Proof} For every object $Y \in \Cal A$,
$\Cal A (Y , \_)$ is the image of $Y$ via the functor
$\hom$, while $< Y^*, \_ >$
is the image of $Y$ via the functor  $\Cal I \circ \Cal E \circ \hom$. But
by Theorem 2.1,  there is a natural equivalence 
$\Cal I \circ \Cal E \simeq id_{\Cal A ^*}$.
Furthermore, these equivalences are natural in $Y$ which proves the
proposition.
\qed
\enddemo

Note that this implies that the inner product is symmetric in the sense 
that 
$$
<Y, Z> \simeq \Cal A ( Y^*, Z) = \Cal A ^{op} (Z,  Y^*) 
	\simeq \Cal A (Z^*,  Y) \simeq < Z, Y>
$$
are naturally isomorphic. 
The first and last isomorphisms are given by Proposition \sect .3 where
we identify $ (Y ^* ) ^* $ with $Y$ and $( Z^*) ^*$ with $Z$ by Proposition
\sect .2.  The middle equation is the usual identification of morphisms
in $\Cal A$ and $\Cal A ^{op}$ as $k$-modules. Using the geometric
definitions, however, we can say  more. Note that $\< \_ , \_\>$ is naturally
equivalent to $tr \circ \fu$ where $tr $ is the functor associated to
the disk $D$ rotated by $180^o$ degrees. The braiding $\sigma$ induces
a natural isomorphism $tr (\sigma ):
< \_ , \_> \simeq <\_ , \_ > \circ \tau$ the
square of which is   $\< \_, \theta ^{-2} \_ \>$ (see Figure 3.4). 
Hence, $tr (\sigma ) \circ  
( 1_{id_\Cal A} \otimes \theta)$ is its own inverse which proves the following.

\vskip .2in
\proclaim{Proposition \sect .4 }
There is a natural equivalence
$
<\_, \_ > \simeq <\_, \_> \circ \tau
$
such that its square  is  the identity.
\qed
\endproclaim 

\vskip 1.5in
\hfil{\it Figure \sect .4: $ tr (\sigma )= \< \_ , \theta ^{-1} \> $.}\hfil

\vskip .2in
\subheading{\sect .3 Adjoints }
Let $\Cal B = \Cal A ( n)$. The inner product on $\Cal A$
extends to an inner product $< \_ , \_ > _\Cal B$ on $\Cal B$ such that
the associated functors $\Cal E $ and $\Cal I$ from the previous section
have compatible extensions to mutually inverse functors on $\Cal B$.
Let $\Cal C = \Cal A ( m)$ and consider the following map 
of categories
$$
ad: [\Cal B, \Cal C] \lra [\Cal C ^*, \Cal B ^*] \overset {(\Cal E, \Cal E )}
\to \lra [\Cal C, \Cal B]
$$
where the first is the dual functor  (A.2) and the second is the equivalence
induced by $\Cal E$ (A.3).
By definition there are natural isomorphisms
$$
<F(X), Y>_\Cal C \simeq < X, ad (F) (Y) > _\Cal B
$$
for all $X\in \Cal B$ and $Y \in \Cal C$. 
As the inner product is non-degenerate,
this defines $ad (F)$ up to natural isomorphism. Let $\bar \Sigma$
denote the cobordism $\Sigma $ rotated by $180^o$ degrees.
Then, from Figure 3.5 we have natural isomorphisms $<F_\Sigma (X), Y>
\simeq < X, F_{\bar \Sigma} (Y)>$, which proves the following.

\vskip .2in
\proclaim{ Lemma \sect .5}
For any cobordism $\Sigma $, there is a  natural isomorphism 
$
ad (F_\Sigma) \simeq F_{\bar \Sigma}.
$
\qed
\endproclaim

\vskip 2.5in
\hfil{\it Figure \sect .5: $ad (F_\Sigma ) \simeq F_{\bar \Sigma }$.}\hfil

\vskip .2in
We may look at the adjoint from a slightly different point of view.
Using A.3,
the involution $^*$ induces a natural equivalence 
$$
(^*, ^*): [\Cal B, \Cal C]
\simeq [\Cal B ^{op}, \Cal C ^{op}].
$$
We denote the image of a functor $F$ by $^* F^*$.
To compare $ ad (F)$ and $^* F ^*$ consider the following natural 
equivalences of categories (see A.4):
$$
\aligned
[\Cal C, \Cal B] &\overset \simeq \to \lra (\Cal B ^{op} \otimes \Cal C) ^*
\quad \text { via } \quad  H \mapsto \Cal B (\_ , H \_ )		\\
[\Cal B ^{op}, \Cal C ^{op}] &\overset \simeq \to \lra (\Cal C \otimes \Cal B
		^{op} )^* \quad \text { via } \quad G\mapsto \Cal C (G\_ ,\_ )
\endaligned
$$
By Proposition 3.3, $H$ and $G$ represent isomorphic functors in $(\Cal B
^{op} \otimes \Cal C) ^*= (\Cal C \otimes \Cal B ^{op} )^*$ if there are
natural equivalences  $<(G \_ ) ^* , \_ > \simeq  <\_ ^* , H \_ >$,
that is $H \simeq ad (^* G ^*)$. For $H= ad (F)$ and $G = ^*F^*$ 
this implies

\vskip .2in
\proclaim {Proposition \sect .6}
$ad (F)$ and $^* F ^*$ are isomorphic as elements in $\Cal B ^{op} \otimes
\Cal C^*$.
\qed
\endproclaim

\vskip .2in
Using the duality described in A.1,
one may also want to compare 
$F \in [\Cal B, \Cal C ] = [\Cal B ^{op}, \Cal C ^{op} ] ^{op} $ with
$^* F ^* \in [\Cal B ^{op} , \Cal C ^{op}]$. 
In general, there is no relation  unless we make further assumptions
(see also Section 5).
In the case when $F= F_\Sigma $ for some endomorphism $ \Sigma 
=\tilde \Sigma \in \Cal S 
(n,n)$, however,  we have canonical
isomorphisms
$$
F_\Sigma = F_{\bar \Sigma } \simeq ad (F_ \Sigma ) 
\simeq ^*F_\Sigma  ^*.
$$
and under these isomorphisms the twist $\theta$ as morphism of the 
idenitity functor $id_\Cal A \simeq F_I$  is mapped to itself.
This is illustrated in Figure \sect .6.

\vskip .2in
\proclaim {Corollary \sect .7}
The twist $\theta$ is self-dual: $\theta _Y ^* = \theta _{Y^*}$.
\qed
\endproclaim

\vskip 1.5in
\hfil {\it Figure \sect .6: $\< \theta , \_ \> = \< \_, \theta \>$.} \hfil

\vskip .2in
Using  similar methods as  above, the following natural equivalences 
can easily be found. Here $\triangle $ is defined to 
be $ ad (\fu) \simeq F_{\bar P} : \Cal A \to \Cal A \hat \otimes \Cal A $. 
$$
\aligned
< \_, \_ \fu \_ > &\simeq < \_ \fu \_, \_> 	\\
<\_, \_> \,  \simeq \,  \tra  \circ \fu \quad &\text {and} \quad 
	<1, \_ >  \, \simeq \, \tra \simeq <\_, 1>	\\
(<\_, \_ > \otimes id_\Cal A ) \circ ( id_\Cal A \otimes \triangle)
 \, \simeq \, 
& \fu  \, \simeq \, 
(id_\Cal A\otimes <\_, \_ >) \circ (\triangle \otimes id_\Cal A)
\endaligned \tag"\sect.8"
$$

\vskip .2in
We remark here that in a monoidal category $\Cal C$ with involution, 
we can take
Proposition \sect .3 as the definition of an inner product:
$< X,Y> := \Cal C (X^* ,Y)$. Simi\-lary, the above formula gives a well
defined trace.
With these definitions, 
the top equation
is equivalent to the condition that there are natural
isomorphisms
$\Cal C(X^*, Y\fu Z) \simeq \Cal C ((X\fu Y)^*, Z)$.
Borrowing a term from algebra, we call  a braided monoidal  category with
these properties and a compatible twist Frobenius.

Conversely, note  that every monoidal category $\Cal C$ with duals satisfies
the top equation in \sect .8 (hence, every tortile category [20]
is Frobenius);
$Y$ and $Y^*$ are right dual to each other if and only if they are
right adjoints in the bicategory (or lax 2-category) with one object 
associated to $\Cal C$. It follows then from well-known results [8; I,6]
that in that case there are natural isomorphisms
$\Cal  C(X^*, Y \fu Z) \simeq \Cal C (Y^*  \fu X^*, Z) \simeq \Cal C ((X\fu 
Y)^*, Z).$

\vskip .4in

\redefine\sect{4}

\subheading{\sect . The Verlinde algebra} 

\vskip .1in
We restrict our attention now to the case when $k$ is an algebraically closed
field. (In most applications $k= \Bbb C$.)  
We show here that the irreducibles
define for us a finite set of \lq labels' with involution. Without any
further assumptions we cannot deduce that $\bold 1$ is irreducible 
or self-dual.

\vskip .2in
\subheading { \sect .1 Simple objects and the dimension of $\Cal A$}
It follows from Theorem 2.5 that $\Cal A $ is an Abelian category.
Recall the notation
$E_A = (\Sigma _{i=1} ^n P_i \otimes Q_i; \pi)$. Up to isomorphism,
we may assume that the $Q_i$ are irreducible and non-isomorphic. Then by
Schur's lemma,

\vskip .2in
\proclaim{Lemma \sect .1 } 
$ \Cal A (Q_i, Q_j) = \delta _{ij} k$
for all $i, j = 1, \dots, n$. 
\qed
\endproclaim

\vskip .1in
As we are working over an algebraically closed field, the irreducible
$A\otimes A$-modules are the tensor products of irreducible $A$-modules.
Hence, $\Cal A \otimes \Cal A$ is already idempotent complete and we
may assume that $E_C = \Sigma _{i=1}^n P_i \otimes Q_i$ with $\pi$
the
identity.

\vskip .2in
\proclaim{Proposition \sect .2 } 
For all $i= 1, \dots, n$, $  Q_i^* \simeq P_i$, and in particular $E^*_C
\simeq E_C$.
\endproclaim

\demo {Proof}
By Theorem 2.1 , $Q_i \simeq \Sigma _j <Q_i, P_j> \otimes Q_j$.
Hence, by the Krull-Schmidt theorem  and above lemma,
$<Q_i, P_j> = \delta _{ij} k$ and thus by Proposition 3.3,
$\Cal A ( Q_i^* , P_j) = \delta _{ij} k$. As each $P_j$ is isomorphic
to a sum of copies of the $Q_l$'s, this forces $P_i \simeq  Q_i^*$.
\qed
\enddemo

\vskip .2in
\proclaim{Corollary \sect .3 } If $O$ denotes a torus, then $E_O \simeq  k^n$,
where $n$ is the number of non-isomorphic irreducibles $Q_i$,
which we call the dimension of
$\Cal A $.
\endproclaim

\demo{Proof}
$E_T = F_T (k) \simeq <\_ , \_> (E_C) = \Sigma _{i=1} ^n < P_i , Q_i> 
\simeq \Sigma _{i=1}^n\, \Cal A ( P^* _i, Q_i ) \simeq \Sigma ^n _{i=1} \,
\Cal A (Q_i, Q_i)
= k ^n$.
\qed
\enddemo

\vskip .2in
We introduce the notation $Q_\alpha$ 
with $\alpha \in  I =\{ 1, \dots , n\}$.
$Q_\alpha ^* \simeq P_\alpha $ is again simple and hence must
be isomorphic to some $Q_{\bar \alpha}$.
By Theorem 2.1 and Proposition 3.3,
every $Y \in \Cal A $ is isomorphic to $\Sigma _\alpha
V_Y^\alpha
\otimes Q_\alpha $ where $V_Y^\alpha = \Cal A ( Q_ \alpha ,Y)$. If $Y= Q_\beta$
we write for short $V_\beta ^\alpha$ and $n^\alpha _\beta$ for its dimension.
Lemma \sect .1
then may be stated as $V_\beta
^\alpha  = \delta _{\alpha \beta} k$.
Note that $\Cal A (Q_\alpha, Q_\alpha)$
has a canonical basis $\{ id _{Q_\alpha} \}$. 
Hence,  by Lemma \sect .1, the the dual space $(V^\alpha _Y )^*$ can 
be identified with $\Cal A ( Y , Q_\alpha)$.

\vskip .2in
\subheading{\sect .2 Verlinde algebra}
We now analyse the multiplicative structure induced by  fusion and
introduce the 
notation
$$
Q_\alpha \fu Q_\beta 
\simeq \Sigma _\gamma \,  V^\gamma _{\alpha \beta} \otimes Q_
		\gamma
$$
where $V^\gamma _{\alpha \beta} = 
\Cal A ( Q_\gamma, Q_\alpha *Q_\beta)$ is 
of dimension $n^\gamma_{\alpha, \beta}$.
The $\Bbb Z$-algebra $\Xi$
with free generators $Q_\alpha, \, \alpha \in I$, and   multiplication
determined by the $n^\gamma_{\alpha, \beta}$ is called the Verlinde algebra.
$\Xi$ is associative and commutative,
$$
n^\gamma_{\alpha , \beta} = n^\gamma _ {\beta, \alpha}
\tag "\sect .1"
$$
as $\Cal A $ is braided monoidal. Furthermore, $\Xi$ has a $\Bbb Z$-linear
involution taking $Q_\alpha$ to $Q _{\bar \alpha}$, and as $\Cal A $ is 
Frobenius,
$$
n^\gamma _{\alpha , \beta} = n^{\bar \beta} _{\bar \gamma , \alpha}.
\tag "\sect .2"
$$
For by Propositions 3.3 and the top equation of 3.8,
we have the following sequence of vector
space isomorphisms:
$$
\aligned
\Cal A ( Q_\gamma, Q_\alpha \fu  Q_\beta)& \simeq <Q_{\bar \gamma},
Q_\alpha \fu  Q _\beta > \simeq < Q_{\bar \gamma} \fu Q_\alpha , Q _\beta>
 \simeq \Cal A ((
Q_ {\bar \gamma}  \fu Q _\alpha)^*, Q _\beta )		\\
&\simeq \Cal A ^{op} 
	(Q_{\bar \gamma} \fu  Q _\alpha,
Q_{\bar \beta}) = \Cal A (Q_{\bar \beta }, Q_{\bar \gamma} \fu  Q _\alpha). 
\endaligned
$$

\vskip .2in
Now let $\bold
1 \simeq  \Sigma _{i=1} ^r \, 
V _{\beta _i} Q_{\beta _i}$ with $ \beta _i \in I$ and $\eta _{\beta _i}
:=\dim V_{\beta _i} >0$.
From the isomorphism $Q_\alpha \fu  \bold 1 \simeq Q_\alpha$, it follows that
$$
\aligned
&\Sigma _i \, \eta _{\beta _i} n^\alpha _{\alpha, {\beta _i}} = 1 
 \\
&\Sigma _i \, \eta _{\beta _i} n^\gamma _{\alpha, \beta _i} = 0  
\quad \text { for }
	\quad \gamma \neq \alpha.	
\endaligned
\tag"\sect .3"
$$
From the second equation, 
we have that $n^\gamma _{\alpha, \beta _i} = 0$
for all $\beta_i$ and $\gamma \neq \alpha$.
Now take $\alpha = \beta _j$. The  first equation then reads $\Sigma _i \,  
\eta _{\beta _i} n^{\beta _j } _{ \beta _j ,  \beta _i} =1$. 
But we just saw that
$n^{\beta _i}
_{\beta _i, \beta _j} = n^{\beta _i} _{\beta _j , \beta _i }= 0$ if
$i\neq j$. Hence,  
$$
n^{\beta _i} _{\beta _i , \beta _i} =1 \quad \text { and }
\quad \bold 1 \simeq  Q _{\beta _1} \oplus  \dots \oplus  Q _{\beta _r}.
$$
For any  $\alpha \in I$, the first equation gives us that
there exist a unique $i = i (\alpha)$
such that $n^\alpha _{\alpha, \beta _{i (\alpha)}} = 1$, and all others are 
zero.
Hence,  define $I_i = \{ \alpha \, | \,  n^\alpha _{\alpha, \beta _i} = 1 \}$
and write $I$ as the union of disjoint sets
$$
I = I_1 \cup \dots \cup I_r.
$$

\vskip .2in
We will now show that the subsets $I_i$ are closed under involution and
multiplication.
Indeed by equations \sect .1 and
\sect .2, $n^\alpha _{\alpha , \beta _i}
= n^{\bar \alpha} _{ \bar \alpha, \beta _i}$,  and therefore $\alpha $ and $\bar
\alpha $ belong to the same subset $I_i$. 
To see that the $I_i$'s are closed under multiplication,
let $\alpha \in I_i$ and $\alpha ' \in I_{i'}$ with $i\neq i'$ and 
evaluate $Q _{\beta _i} \fu  ( Q _\alpha \fu  Q _ {\alpha '}) \simeq 
Q _{\alpha} \fu  Q _{\alpha '} \simeq ( Q _{\alpha} \fu Q _{\alpha '}) 
\fu  Q _{\beta
_{i'}}$ in two different ways. The former is a sum of
$Q_\gamma$'s  with $\gamma \in I _i$ and the latter is a sum of
$Q _\gamma$'s with $\gamma \in I _{i'}$. As $I_i$ and $I_{i'}$ are
disjoint, this implies that the expression is zero and
in particular that  $n^\gamma _{\alpha , \alpha '} = 0$.

\vskip .2in
\proclaim{ Proposition \sect .4}
The Verlinde algebra $\Xi$
is a direct product of $r= \dim \Cal A (\bold 1, \bold 1)$
unital, finitely generated   
algebras  $\Xi _i = \Bbb Z [Q_\alpha \, | \, \alpha \in I _i]$ with involution. 
\qed
\endproclaim

\vskip .1in
Note that we do not claim that the involution is an algebra homomorphism
or that it  fixes the unit (see also Remark 5.2).

\vskip .4in
\redefine\sect{ 5}

\vskip .2in
\subheading{\sect . Orientation and  modular functors}

\vskip .2in
Most definitions of a topological field theory or a modular functor 
contain a condition on how the data transforms when the orientation of 
a surface is reversed.
So far, it does not make sense in our theory to reverse orientation as all
surfaces in $\Cal S$
are 
oriented such that the normal vector is outward pointing. It is natural however
to consider the following canonical extension.

\vskip.2in
\subheading{\sect .1. Extension to the oriented category}
Enlarge the category $\Cal S $ to $\Cal S ^c$ so that its objects
are pairs $(n, w)$ where $w$ is one of the $2^n$ possible orientations
of $n$.
Similarly, cobordisms
come now with a chosen orientation (determined by whether the normal
vector is outward or inward pointing). The orientation of the boundary
is the induced orientation. $\Cal S$ can be identified with 
a subcategory of $\Cal S ^c$.
$\Cal S^c$ has a slightly more
complex structure than $\Cal S$ in that it is most natural to think
of it as a double category [12] which has two types of 1-morphisms. 
In addition to the cobordisms, which we think of the horizontal 1-morphisms,
we also have vertical 1-morphisms induced  by reflection $c$ in the $x$-axis
which maps a
circle to itself with opposite orientation.
\footnote{
Indeed, one might want to impose a double category structure on $\Cal S $
by including all orientation preserving diffeomorphisms of the 1-manifolds
as vertical 1-morphims.
But  up to homotopy the diffeomorphism group of $n$ copies of a
circle is just the wreath product $\Sigma _n \wr S^1$ of the symmetric
group on $n$ letters and the circle.  This is already represented
by crossed   cylinders and their twists,  the 2-morphisms $\theta$.
Note that the relation between the 2-morphisms and 1-morphisms is
given by taking the classifying spaces. Here we have $B \Bbb Z \simeq
S^1$. By Proposition 3.6, we can  represent  conjugation
$c$ by the bend cylinder $C$. }
The 2-morphisms are not only those orientation preserving
diffeomorphisms which fix the boundary (2-cells with identity 
vertical maps) but also those which restrict to given maps
on the boundary. These give 2-cells of the form:
$$
\CD
(n, w_1)	@> F_{\Sigma _1} >>	(m, w_2)	\\
@V \bar  c_1 VV				@V \bar c _2 VV	\\
(n, \bar w_1)	@> F_{\Sigma _2}>>	(m, \bar w_2).
\endCD
$$
When the surfaces are connected
any 2-cell with non-identity vertical maps 
can uniquely  be
written  as a diffeomorphism  that fixes the boundary 
composed with  the reflection  $\gamma$ 
in the $xz$-plane. Clearly, the monoidal 
structure of $\Cal S $ extends to a monoidal structure of $\Cal S ^c$.
We denote 
a cobordism $\Sigma$ with the opposite orientation by $\bar \Sigma$.
\footnote{
The following extension of $\Cal A$ 
justifies the use of the same symbol for a surface rotated by $180^o$ and
its conjugate
in the light of  Lemma 3.5 and Proposition 3.6.} 

\vskip .2in
Define now the following natural extension of $\Cal A $ to a 
monoidal functor $\Cal A^c : \Cal S ^c \to \wA $ by setting
$$
\aligned
\Cal A ^c (\bar S^1) &:= \Cal A (S^1) ^{op} = \Cal A ^{op}	\\
\Cal A ^c(c) &:=  ^*		\\
\Cal A ^c( \bar \Sigma )&:= F_{\bar \Sigma } = ^*  F_\Sigma ^*. 
\endaligned
$$ 
With these definitions, the  2-cell
$\gamma$ corresponds to the obvious   natural equivalence $^* \circ F_\Sigma
\simeq ^* \circ F_\Sigma \circ ^* \circ ^* = F_{\bar \Sigma } \circ  ^*$.

\vskip .2in
\proclaim{ Lemma \sect .1}
For closed surfaes, $  E _{\bar \Sigma }$ is  canonically
isomorphic to $E_\Sigma ^*$.
\endproclaim

\demo{Proof}
The canonical object  $k \in \hat k$ 
is canonically isomorphic to its dual in $\hat k ^*$.
Hence, there is a canonical isomorphism
$E_\Sigma ^* = F_\Sigma ^* (k)  \simeq
F_\Sigma ^* (k^*) = F _{\bar \Sigma } (k) = E_{\bar \Sigma }$.
\qed
\enddemo

\vskip .2in
\demo{Remark \sect .2}
As mentioned in Section 3, in general there is no relation between
$F_\Sigma$ and $F_{\bar \Sigma}$ so that $F_{\bar P}$ defines another
braided monoidal structure $\bar \fu$
on $\Cal A$ with unit $\bold 1 ^* = E_{\bar D}$. 
By Proposition 3.7, 
$\theta$ is a compatible twist.
The two monoidal structures are related by natural isomorphisms
$(Y\fu Z)^*  \simeq Z^* \bar \fu Y^*$.  Therefore, unless
$\fu $ and $\bar \fu $ are naturally equivalent, $\Cal A$ cannot be 
tortile in the sense of [20]. 
Conversely, if $(Y\fu Z)^* \simeq Z^*\fu Y^*$ via natural isomorphisms
then $\Cal A$ is tortile by the results in Section 3 [10].
In this case, note that $F_\Sigma \simeq F_{\bar \Sigma}$ for any
surface $\Sigma$ and by Lemma \sect .1 every $E_\Sigma$ is canonically
isomorphic to its dual, i.e. is an inner product space.
When $k$ is algebraically closed, $\Cal A $ would then
be   an artinian, semi-simple
tortile category in which the morphism space of the irreducibles is 
1-dimensional.
Yetter [24] constructs explicitly a state sum 3-manifold invariant from
such $\Cal A$  which generalizes that of
Tuarev and Viro [22].

\enddemo

\vskip .2in
\subheading{\sect .2. Modular Functor}
We briefly indicate here how one can construct 
a modular functor in the sense of [19] and [14] from $\Cal A$ when $k= \Bbb C$.
Modular functors are of importance because of their one to one correspondence
with 2+1 dimensional field theories [14] [23].

\vskip .2in
In the previous section we saw that the set of irreducible objects
up to isomorphisms is a finite set $I$ with involution. 
By 
Proposition 4.4, we may restrict our attention for the moment
to the case when
$r= \dim \Cal A (\bold 1, \bold 1) =1$, in other words, when
$\bold 1 = Q_0$ is irreducible. We make the following
{\it non-triviality assumption}:
$$
E_{S^2} \simeq <E_D, E_D> \simeq \Cal A (\bold 1^*, \bold 1) \neq 0.
$$
As $\bold 1^* \simeq Q_{\bar 0}$ is irreducible, this then implies 
that $E_{S^2} \simeq \Bbb C$ and $\bold 1 ^* \simeq \bold 1$, that is
$\bar 0 = 0$ for the vacuum element $0 \in I$.

\vskip .2in
Given a surface $\Sigma $ with boundary $n$ and a colouring  $\sigma : n
\to I$, define $V(\Sigma, \sigma) := F_\Sigma (Q_\sigma)$ where $Q_
\sigma = Q_{\sigma (1)}
\otimes \dots \otimes Q _{\sigma(n)}$. Then the mapping class group
$\Gamma (\Sigma)$ acts on $V(\Sigma , \sigma)$ as it acts on $F_\Sigma$,
and by definition of $F_
{\bar \Sigma}$, $\gamma$ defines a natural isomorphism
$V(\Sigma , \sigma )^* \simeq V(\bar \Sigma , \bar \sigma)$.
Also, by our non-triviality assumption and equations 4.3, 
$V(D, \sigma) $ is 0 unless $\sigma = 0$ in which case it is $\Bbb C$.
Because the functor $\Cal A$ is a monoidal functor, 
if $\Sigma = \Sigma _1 \cup \Sigma _2$ and $\sigma =
\sigma _  1\cup \sigma _2$ then $V(\Sigma , \sigma ) \simeq
V(\Sigma _1, \sigma _1) \otimes V (\Sigma _2, \sigma _2)$.
The gluing law is equally implied; let $\Sigma ' $ be $\Sigma $ 
with the boundary components $n-1$ and $n$ identified.
This correspondes to gluing a cylinder onto these boundary components.
Then $V( \Sigma ' , \sigma ) \simeq  F_\Sigma ( Q_\sigma \otimes E_C)
\simeq \bigoplus _{\alpha \in I} F_\Sigma  (Q_\sigma \otimes Q_{\bar \alpha}
\otimes Q_\alpha) = \bigoplus _{\alpha \in I}
V(\Sigma , \sigma \cup \{\bar \alpha , \alpha \})$.

\vskip .2in
Finally, the twist $\theta$ defines for each label $\alpha \in I$
a non-zero complex number $h_\alpha :=
\theta _{Q_\alpha} \in \Cal A  (Q_\alpha,
Q_\alpha ) = \Bbb C$. 
As $\theta $ is self-dual by Proposition 3.7, 
$h_{\bar \alpha}  = h_\alpha$.  Also,
$h_0 = \theta _{E_D} = 1$ as $\Cal A $ is balanced by Proposition 3.1.
In summary, we have proved:

\vskip .2in
\proclaim{Theorem \sect .2} 
Every $\Bbb C$-linear category $\Cal A $ with $\Cal S$-structure
determines $\dim E_{S^2}$ many   modular functors.
\qed
\endproclaim
\redefine\sect{A}

\vskip .4in
\subheading {Appendix: Linear categories }

\vskip .1in
We review some essential notions from the theory of linear categories. 
Let $k$ be a commutative ring.  A $k$-linear category (or a $k$-category for
short) is a category equipped with a $k$-module structure on each morphism
set such that composition is $k$-bilinear. A $k$-functor between
two $k$-categories $\Cal A $ and $\Cal A'$ is a functor which 
is $k$-linear on the morphisms, i.e. it maps $\Cal A (p,q)$ linearly
to $\Cal A' (Fp, Fq)$. 
When both $\Cal A $ and $\Cal A'$ are small,
we denote the category of $k$-linear functors from
$\Cal A $ to $\Cal A'$ by $[\Cal A , \Cal A']$ or $(\Cal A ' ) ^{\Cal A }$.

\vskip .2in
\demo{Tensor Product}
Given two $k$-categories $\Cal A$ and $\Cal A'$, we can form their (naive)
tensor product $\Cal A \otimes \Cal A'$ over $k$ by setting
$$
\aligned
&\text {ob} ( \Cal A \otimes \Cal A') = \text {ob} (\Cal A) \times \text{ob} 
	(\Cal A')	\\
&(\Cal A \otimes \Cal A') ((p,p'), (q,q'))= \Cal A (p,q) \otimes \Cal A'(p', q')
\endaligned
$$
with composition defined by $(f\otimes f') (g\otimes g') = (fg \otimes f'g')$.
Clearly, $\Cal A \otimes \Cal A'$ is again a $k$-category. We have the following
canonical isomorphisms of categories
$$
\aligned
&(\Cal A \otimes \Cal A') 
	\otimes \Cal A'' \simeq \Cal A \otimes (\Cal A' \otimes
	\Cal A'' )	\\ 
&\Cal A \otimes \Cal A' \simeq \Cal A' \otimes \Cal A	\\
&\Cal A \otimes k \simeq \Cal A \simeq k \otimes \Cal A
\endaligned
$$
where $k$ denotes the category with one object and morphism set $k$.

\vskip .1in
We will always consider small $k$-categories also called $k$-algebroids. These
have been studied extensively, see for example [15] and [16], and may be
thought of
as algebras with several objects.  
The category $\A$ of all small $k$-categories, $k$-linear functors and 
natural transformations is a strict 2-category, and the tensor product
induces a symmetric monoidal structure on the underlying category. Indeed it 
is easily checked that $\A$ is a monoidal 2-category in the sense of 
[7; 2.6].
\enddemo

\vskip .2in
\demo{Additive and idempotent completion}
The categories of interest to us should have the additional properties that any
two objects have a coproduct and that idempotents split. In other words, the
categories should be additive and idempotent complete.  Such categories
are called amenable or also Karoubian.

\vskip .1in
Given any $k$-category $\Cal A$ there is  an  amenable completion $\hat
\Cal  A$ which is unique
up to equivalence relative to  $\Cal A$. More concretely, an 
additive completion of $\Cal A$ is given by the category $Mat (\Cal A)$ the
objects of which are finite sequences of objects in $\Cal A$ and the morphisms
of which are matrices of morphisms in $\Cal A$.  An idempotent
completion of $\Cal A$ can be constructed as follows. The objects of the new 
category are the idempotents $(p,e)$ with $e^2 =e \in \Cal A (p,p)$,
and a morphism between two
idempotents  is a triple $(e',f,e)$ where $f$ is a morphism of
$ \Cal A$ satisfying the equation $e'f = f =  fe$.  Composition is defined 
by 
$$
(e'',f',e') (e',f,e)= (e'',f'f,e)
$$
Thus
$\widehat {Mat (\Cal A)}$ is an amenable completion of $\Cal A$. For example, 
consider a $k$-algebra $A$
as a category $\Cal A$ with one object,  then its amenable completion
$\widehat {Mat (\Cal A)}$ is 
equivalent to the category of finitely generated projective $A$-modules.

\vskip .1in
Every functor  from
a $k$-category to an amenable category 
extends to its
completion
and any natural transformation of two such functors extends uniquely to a
natural transformation of the extended functors. In other words, the
categories $[ \Cal A , \Cal A']$ and $[\hat {\Cal A}, {\Cal A'}]$ are
naturally equivalent when $\Cal A '$ is an amenable category.
Thus the category $\Cal A$ is Morita equivalent to its
completion
in the following sense.  

\vskip .1in
Let $Mod (k)$ be the category of $k$-modules, and
define
$$
Mod (\Cal A) := Mod (k) ^{\Cal A} = [\Cal A, Mod (k)]
$$
to be the category of $k$-linear functors from $\Cal A $ to $Mod (k)$.
Then two small $k$-categories $\Cal A$ and $\Cal B$ are Morita equivalent
if their module categories $Mod (\Cal A)$ and $Mod (\Cal B)$ are 
equivalent.
If $p\in \Cal A$ is a fixed object, then
$$
\Cal A(p,\, \,) : \Cal A \lra Mod( k) \, \text { via } \, q \lra \Cal A (p,q)
$$
defines a natural element in $Mod (\Cal A )$. Letting $p$ vary, we get a 
contravariant functor 
$$
\hom : \Cal A \lra Mod(\Cal A ).
$$
$Mod (\Cal A)$ should not be confused with the smaller category $\Cal A ^*
:= [\Cal A , \hat k]$ which we call the dual of the category $\Cal A$.
The hom-functor  factors through $\Cal A ^*$ if and only if the morphism
sets $\Cal A (p,q)$ are finitely generated projective $k$-modules for all
$p$ and $q$.

\vskip .1in
Now, let  $\wA  \subset \A$ denote  the full subcategory of 
amenable $k$-categories. The as\-sign\-ment 
$\Cal A \to \widehat {Mat (\Cal A )}$
induces a  2-functor $\A \to \wA $. 
Replacing the tensor product
by the completed tensor product
$$
\Cal A \hat \otimes \Cal A' := \widehat {(\Cal A \otimes \Cal A')}
$$
$\wA$ is again a monoidal 2-category. Its unit is $\hat k$, the category of
finitely generated projective $k$-modules.
An object of $\Cal A \hat \otimes \Cal A'$ may
be written as $(\Sigma _{i=1}^n (p_i, q_i); e)$ where 
$$
e^2=e \in
\Cal A  \otimes \Cal A' (\Sigma _i (p_i,q_i), \Sigma _i (p_i,q_i))
	= \Sigma _{i,j} \Cal A (p_i, p_j) \otimes \Cal A' (q_i,q_j).
$$
\enddemo

\vskip .2in
\demo{Involution and dual}
We list here some elementary facts
that clarify the functoriality  of various maps
between functor categories.
\enddemo

\vskip .1in
\demo{\sect .1}
If $F:\Cal A \to \Cal B$ is a covariant functor, define $F^{op}: \Cal A ^{op}
\to \Cal B ^{op}$ that acts on objects like $F$ and takes $f^{op}:p \to q$
to $F(f)^{op}:F(p) \to F(q) $.
The assignment $F \to F^{op}$ defines a contravariant functor
$$
{\,}^{op}: [\Cal A, \Cal B] \lra [\Cal A ^{op} , \Cal B ^{op}].
$$
Indeed, a natural transformation $\phi: F_1 \to F_2$ is mapped to
$(\phi ^{op})_p := (\phi_p)^{op} : F_2 (p) \to F_1 (p)$.
Clearly, this is an involution in the sense that ${\,}^{op} \circ {\,}^{op}
= id _{[\Cal A, \Cal B]}$.
\enddemo

\vskip .1in
\demo{\sect .2}
There is a naturally defined, covariant  dual functor
$$
*: [\Cal A , \Cal B] \lra [Mod (\Cal B), Mod (\Cal A)]
$$
taking a functor $F$ to $F^*(G) := G\circ F$ where $G: \Cal B \to Mod (k)$
is an element in $Mod (\Cal B)$. Clearly, a natural transformation 
$\phi : F_1 \to F_2$ is just send to $\phi ^* : F_1 ^* \to F_2 ^*$
where
$\phi^* _G (p):= G ( \phi _p)$.
\enddemo

\vskip .1in
\demo {\sect .3}
Let $H:\Cal A \to \Cal A'$ and $ G: \Cal B \to \Cal B '$ be covariant
equivalences with $H\circ H' \simeq id _{\Cal A '}$ and $G' \circ G
\simeq id _{\Cal B}$.  Then 
$$
(H,G): [\Cal A , \Cal B] \lra [\Cal A ' , \Cal B']
$$
defined by $F \to G \circ F \circ H'$ is a covariant equivalence
of categories.   
Indeed, a natural transformation $\phi: F_1 \to F_2$ is taken to
the natural transformation $(G\phi H')_p = G ( \phi _{ H'(p)})$.
\enddemo

\vskip .1in
\demo {\sect .4}
If $\Cal C$ is an additive $k$-category  then there are natural 
isomorphisms of $k$-categories
$$
(\Cal C ^\Cal A ) ^\Cal B \simeq \Cal C ^{\Cal A \otimes \Cal B }.
$$
In particular, $Mod(\Cal A ) ^\Cal B \simeq Mod ( \Cal A \otimes \Cal B)$
and $[\Cal A , \Cal B ^*] \simeq (\Cal B \otimes \Cal A ) ^*$.
\enddemo

\vskip .4in
\subheading{References}

\parindent =0in

\vskip .2in

[1] M.F. Atiyah,  Topological quantum field theory, Publ. Math. IHES
(Paris)  68 (1989), 175--186.

\vskip .1in
[2] C. Blanchet, N. Habegger, G. Masbaum , P. Vogel,  Topological 
quantum field theories derived from the Kauffman bracket, Topology 34 (1995),
883--929.

\vskip .1in
[3] S. Carmody, \lq \lq Cobordism Categories", D.Phil.-Thesis, Cambridge
University, 1995.

\vskip .1in
[4] R. Cohen, J.D.S. Jones, G. Segal,  Floer homotopy type,
in Floer Memorial Volume, ed. E. Zehnder, Birkh\"auser 1994.

\vskip .1in
[5] F. DeMeyer, E. Ingraham,  Separable algebras over commutative rings,
Springer Verlag LNM  181, (1971).

\vskip .1in
[6] E. Getzler,  Batalin-Vilkovisky Algebras and two-dimensional
topological field theories, Comm. Math. Phys.  159 (1994),
265--285.

\vskip .1in
[7] R. Gorden, A.J. Power, R. Street, Coherence for tricategories,
AMS Memoirs 558 (1995).

\vskip .1in
[8] J.W. Gray, Formal Category Theory: Adjointness for 2-Categories,
Springer LNM  391 (1974).

\vskip .1in
[9] D. Freed, Higher algebraic structures and quantization,
Comm. Math. Phys.  159  (1994), 343--398.

\vskip .1in
[10] M. Hyland, private communications.

\vskip .1in
[11] A. Joyal, R. Street,  Braided tensor categories, Adv. in Math.
102 (1993), 20--78.

\vskip .1in
[12] G.M. Kelly, R. Street, Review of the elements of 2-categories,
Springer Verlag LNM  420 (1974), 75--103.

\vskip .1in
[13] M.M. Kapranov, V.A. Voevodsky, 
2-catgories and Zamolodchikov Tetrahedra
Equations, AMS Proc. 56 (1994), 177--259.

\vskip .1in
[14] M. Kontsevich,
Rational conformal field theories and invariants of 3-di\-men\-sional
manifolds,  Marseille preprint CPT-88/P.2189.

\vskip .1in
[15] B. Mitchell,  Rings with several objects, Adv. in Math.  8 
(1972), 1-161.

\vskip .1in
[16] B. Mitchell,  Separable algebroids, AMS  Memoirs  333 (1985). 

\vskip .1in
[17] N. Y. Reshetikhin, V.G. Turaev,  Invariants of 3-manifolds via link
polynomials and quantum groups, Invent. Math.  103 (1991), 547-97.

\vskip .1in
[18] G. Segal,  
The definition of  a conformal field theory, in preparation.

\vskip .1in
[19] G. Segal,  Conformal field theories and modular functors, in
Proc. Int. Congress of Appl. Math., Swansea 1988.

\vskip .1in
[20] M.C. Shum, Tortile tensor categories, J. Pure Appl. Alg.  93
(1994), 57--110.

\vskip .1in
[21] U. Tillmann,  Discrete models for the category 
of Riemann surfaces, to appear in Math. Proc. Cam. Phil. Soc..

\vskip .1in
[22] V.G. Turaev, O.Y. Viro,  State sum invariants of 3-manifolds
and quantum 6j-symbols, Topology  31 (1992), 865--902.

\vskip .1in
[23] K. Walker,  On Witten's 3-manifold invariants, preprint, 1991.

\vskip .1in
[24] D.N. Yetter,  State sum invariants of 3-manifolds associated
to artinian semisimple tortile categories, Topology and its Appl.  58
(1994), 47--80.

\vskip .3in
\noindent
Mathematical Institute 

\noindent
24-29 St Giles Str.

\noindent
Oxford, OX1 3LB

\noindent
UK

\vskip .05in
\noindent
tillmann\@maths.ox.ac.uk
\enddocument